\documentclass[11pt]{article}
\usepackage{epsf, epsfig, color}
\usepackage{amssymb,latexsym}
\usepackage{mathrsfs}
\usepackage{amsfonts}

\newtheorem{lem}{Lemma}
\newtheorem{theo}{Theorem}

\newtheorem{cor}{Corollary}

\newtheorem{prob}{Problem}

\renewcommand{\theenumi}{\rm (\roman{enumi})}

\newcommand{\proof}
{{\noindent {\em Proof}.\quad}\setcounter{countclaim}{0}
\setcounter{countcase}{0}}
\newcommand{\proofend}{{\hfill$\Box$}}

\newcounter{countfig}

\newcounter{countclaim}
\def\inclaim{\addtocounter{countclaim}{1}
{\noindent {\bf Claim \thecountclaim}: }}

\newcounter{countcase}

\def\setz{{\cal Z}}
\def \setr{{\cal R}}

\def \setgr {\mathscr{R}}

\newcommand{\beeq}{\begin{equation}}
\newcommand{\eneq}{\end{equation}}

\newcommand{\beeqn}{\begin{eqnarray*}}
\newcommand{\eneqn}{\end{eqnarray*}}

\newcommand \nroot[1] 
{{\lceil \frac{27 #1}{11} -\frac{27}{22}\rceil+2\mu(6-#1)}}

\newcommand {\relabel}[1] {\label{#1} \red{[*: #1]}}\newcommand {\rebibitem}[1] {\bibitem{#1} \red{[*: #1]}}
\def\relabel {\label} \def\rebibitem {\bibitem}  

\begin{document}
  
\newcommand{\resection}[1]
{\section{#1}\setcounter{equation}{0}}

\renewcommand{\theequation}{\thesection.\arabic{equation}}

\renewcommand{\labelenumi}{\rm(\roman{enumi})}

\baselineskip 0.53 cm

\title {On graphs whose flow polynomials have real roots only}

\author
{ 
Fengming Dong\thanks{This paper was partially supported 
by NTU AcRF project (RP 3/16 DFM) of Singapore.}\\
\small Mathematics and Mathematics Education\\
\small National Institute of Education\\
\small Nanyang Technological University, Singapore
}

\date{}

\maketitle

\begin{abstract}
Let $G$ be a bridgeless graph. 
In 2011 Kung and Royle showed that  
all roots of the flow polynomial $F(G,\lambda)$ of $G$
are integers 
if and only if $G$ is the dual of a chordal and plane graph. 
In this article, we study whether
a bridgeless graph $G$ for which 
$F(G,\lambda)$ has real roots only 
must be the  dual of some chordal and plane graph.
We conclude that the answer of this problem 
for $G$ is positive 
if and only if  $F(G,\lambda)$ does not have any real root in the interval $(1,2)$.
We also prove that 
for any non-separable and $3$-edge connected $G$, 
if $G-e$ is also non-separable for each edge $e$ in $G$
and every $3$-edge-cut of $G$ consists of edges incident with 
some vertex of $G$, 
then all roots of $P(G,\lambda)$ are real if and only if 
either $G\in \{L,Z_3,K_4\}$ 
or $F(G,\lambda)$ contains at least $9$ real roots 
in the interval $(1,2)$,
where $L$ is the graph with one vertex and one loop and 
$Z_3$ is the graph with two vertices and three parallel 
edges joining these two vertices. 
\end{abstract}

\def \MSC {Mathematics Subject Classifications}

{\bf \MSC}: {05C21, 05C31}


\section{Introduction}

The graphs considered in this paper 
are undirected and finite,
and may have loops and parallel edges. 
For any graph $G$, let $V(G), E(G), P(G,\lambda)$ 
and $F(G,\lambda)$ be the set of vertices,
the set of edges, the chromatic polynomial 
and the flow polynomial of $G$.
The roots of $P(G,\lambda)$ and $F(G,\lambda)$ 
are called {\it the chromatic roots}
and {\it the flow roots} of $G$ respectively. 
As $P(G,\lambda)=0$ (resp. $F(G,\lambda)=0$) 
whenever $G$ contains loops (resp. bridges), 
we will assume that $G$ is 
loopless (resp. bridgeless) 
when $P(G,\lambda)$  (resp. $F(G,\lambda)$) is considered.

The chromatic polynomial  $P(G,\lambda)$ of $G$ 
is a function which counts the number of proper $\lambda$-colourings whenever $\lambda$ is a positive integer. 
A chordal graph $G$ is a graph in which every subgraph 
of $G$ induced by a subset of $V(G)$ is not 
isomorphic to any cycle of length larger than 3. 
It is known that if $G$ is chordal, then 
all chromatic roots of $G$ are non-negative integers
(see \cite{dong2005, read1988, rea}). 
Some non-chordal graphs also have this property
(see \cite{dmi1980,dong2005, dong3, dong4, read1975}). 
Meanwhile, there are graphs which have real 
chromatic roots only but also have 
non-integral chromatic roots.
For example, when $s\ge 7$, 
the graph $H_s$ obtained from $K_s$ by 
subdividing a particular edge once is such a graph, as   
\beeq \relabel{eq1-1}
P(H_s,\lambda)=\lambda (\lambda-1)\cdots (\lambda -s+2)
(\lambda^2-s\lambda+2s-3).
\eneq 
However, it is still unknown if there is a 
planar graph $G$ with this property, i.e., 
{\it $G$ has real chromatic roots only 
but also contains non-integral chromatic roots}. 
Due to Tutte~\cite{tut3}, $P(G,\lambda)=\lambda F(G^*,\lambda)$
holds for any connected plane graph $G$, where 
$G^*$ is the dual of $G$.
Thus, equivalently, it is unknown if there is a 
planar graph $G$ which has real flow roots only 
but also has non-integral flow roots.
Actually it is also unknown if there is a non-planar graph 
with this property. 
It is natural to consider the following problem. 

\begin{prob}
\relabel{prob1} 
Is there a bridgeless graph which has real flow roots only 
but also contains non-integral flow roots?
\end{prob}

By the following result due to Kung and Royle \cite{kung},
Problem~\ref{prob1} is equivalent to whether 
there exists a graph $G$ which is not the dual of any plane 
and chordal graph but has real flow roots only. 
If there does not exist any graph asked in Problem~\ref{prob1}, 
then every graph with real flow roots only must be the dual of some chordal and plane graph.

\begin{theo}
[\cite{kung}]
\relabel{kungroyle-theo}
If $G$ is a bridgeless graph, then 
its flow roots are integral
if and only if $G$ is the dual of a chordal and plane graph. 
\end{theo}

In this paper, 
let $\setgr$ be the family of bridgeless  
graphs which have real flow roots only. 
We will focus on graphs in $\setgr$ and 
mainly show that 
for any graph $G\in \setgr$,
all flow roots of $G$ are integers 
if and only if 
$G$ does not contain any real flow roots in the interval 
$(1,2)$.

A vertex $x$ in a connected $G$ is called a 
{\it cut-vertex} if $G-x$ has more components than $G$ has,
where $G-x$ is the graph obtained from $G$ by deleting 
$x$ and all edges incident with $x$. 
A graph $G=(V,E)$ is said to be 
{\it non-separable} if 
either $|E|=|V|=1$ or 
$G$ is connected without loops or cut-vertices.
An edge-cut $S$ of a graph $G=(V,E)$ is the set of edges 
joining verteces in $V_1$ to vertices in $V_2$
for some  partition $\{V_1, V_2\}$ of $V$. 
An edge-cut $S$ of $G$ is said 
to be {\it proper} if $G-S$ has no isolated vertices.

The definition of the flow polynomial of a graph $G$ 
is given in (\ref{eq1-2}).
By the second equality in (\ref{eq1-2}),
$F(G,\lambda)=0$ holds if $G$ contains bridges. 
By the second and the fifth equalities in (\ref{eq1-2}),
$F(G,\lambda)=F(G/e,\lambda)$ if $e$ is one edge 
in a $2$-edge-cut of $G$. 
For this reason, 
the study of flow polynomials can be restricted to 
$3$-edge connected graphs. 
By Lemmas~\ref{block-factor},~\ref{v-edge} and~\ref{2-edge}, 
the flow polynomial of any graph can be expressed as the product 
of flow polynomials of graphs $G$ satisfying the following 
conditions,  divided by $(\lambda-1)^a(\lambda-2)^b$ for some 
non-negative integers $a$ or $b$:
\begin{enumerate}
\item $G$ is non-separable and 
$3$-edge connected;
\item 
$G$ does not contain any proper $3$-edge-cut; and 
\item $G-e$ is non-separable for each edge $e$ in $G$.
\end{enumerate}

Let $\setgr_0$ be the family of those graphs in $\setgr$  
which satisfying conditions (i), (ii) and (iii) above.
By Lemmas~\ref{block-factor},~\ref{v-edge} and~\ref{2-edge},
there exists a graph asked in Problem~\ref{prob1}
belonging to $\setgr$ 
if and only if there exists a graph 
asked in Problem~\ref{prob1} belonging to $\setgr_0$.
Thus the study of Problem~\ref{prob1} can be focused on 
graphs in $\setgr_0$. 

Let $W(G)$ be the set of vertices in a graph $G$ 
of degrees larger than $3$
and let $\bar d(G)$ be the mean 
of degrees of vertices in $W(G)$.
Let $L$ denote the graph with one vertex and one loop
and 
let $Z_k$ denote the graph with two vertices 
and $k$ parallel edges joining these two vertices. 
Our main result in this paper is the following one.

\begin{theo}\relabel{main-th}
Assume that $G=(V,E)$ is any graph in $\setgr$.
\begin{enumerate}
\item
If some flow roots of $G$ are not in the set $\{1,2,3\}$, 
then $|E|\ge |V|+17$ and 
$G$ has at least $9$ flow roots in the interval $(1,2)$.
\item If $G\in \setgr_0$,
then either $G\in \{L, Z_3,K_4\}$
or $G$ has the following properties: 
\begin{enumerate}
\renewcommand{\theenumi}{\roman{enumi}}
\renewcommand{\theenumii}{\rm \theenumi.\arabic{enumii}}
\item\relabel{no1} 
$3\le |W(G)|< \frac {11}{27}|V|+\frac 5{27}$; 
\item\relabel{no2} $G$ contains at least 
$\nroot{|W(G)|}\ge 9$
flow roots in $(1,2)$,
where $\mu(x)$ is the function defined by $\mu(x)=1$ when $x>0$ 
and $\mu(x)=0$ otherwise;

\item\relabel{no3} $\bar d(G)> 14.656-11.656/|W(G)|>10.770$;
\item\relabel{no4} $|V|+8|W(G)|-7\le |E|<(32|V|-49)/5$. 
\end{enumerate}
\end{enumerate}
\end{theo}

Note that 
Theorem~\ref{main-th} (i) is from 
Theorem~\ref{main-lem1-cor}, 
Theorem~\ref{main-th} (ii.1)
is from Theorem~\ref{main-lem1} \ref{cor1-no01} and
Lemma~\ref{main-lem30} \ref{cor1-no1},
Theorem~\ref{main-th} (ii.2) 
are in Theorem~\ref{main-lem1} \ref{cor1-no2},
while Theorem~\ref{main-th} (ii.3) and (ii.4) are from 
Lemma~\ref{main-lem30} \ref{lem30-no3}, ~\ref{lem30-no4}
and~\ref{lem30-no6}.

\noindent {\bf Remark}: Theorem~\ref{main-th} (ii)
implies that Theorem~\ref{kungroyle-theo} holds 
for all graphs in $\setgr_0$.

Interestingly, 
Theorems~\ref{kungroyle-theo} and~\ref{main-th}
imply three equivalent statements on 
a bridgeless graph which has real flow roots only.

\begin{cor}\relabel{sect1-cor}
Let $G\in \setgr$. 
Then the following statements are equivalent:
\begin{enumerate}
\item $G$ is the dual of some chordal and plane graph;
\item  $G$ does not have any flow root in the interval $(1,2)$;
\item each flow root of $G$ is in the set $\{1,2,3\}$.
\end{enumerate}
\end{cor}

\resection{Basic results on flow polynomials}
\relabel{sec2}

Let $G=(V,E)$ be a finite graph with vertex set 
$V$ and edge set $E$ and let $D$ be an orientation of $G$. 
For any finite additive Abelian group $\Gamma$, 
a {\it $\Gamma$-flow} on $D$ is a mapping 
$\phi: E\rightarrow \Gamma$ such that 
\begin{equation}\relabel{sec1-eq1}
\sum_{e\in A^+(v)}\phi(e)=\sum_{e\in A^-(v)}\phi(e)
\end{equation}
holds for every vertex $v$ in $G$,
where 
$A^+(v)$ (resp. $A^-(v)$) is the set of loopless arcs in $D$ 
with tail $v$
(resp. with head $v$).
If $\phi(e)\ne 0$ for all $e\in E$, then 
a $\Gamma$-flow $\phi$ on $D$ is called a {\it nowhere-zero}
$\Gamma$-flow on $D$.
For any integer $q\ge 2$, 
a {\it nowhere-zero $q$-flow} of $G$ is defined to be 
a nowhere-zero ${\mathbb Z}$-flow $\psi$
such that $|\psi(e)|\le q-1$ for all $e\in E$,
where ${\mathbb Z}$ is the additive group consisting of 
all integers. 
Tutte~\cite{tut2} showed 
that $G$ has a nowhere-zero $q$-flow
if and only if it has a nowhere-zero $\Gamma$-flow,
where $q$ is the order of $\Gamma$.

The {\it flow polynomial} $F(G,\lambda)$ of a graph $G$
is a function in $\lambda$ which counts 
the number of nowhere-zero $\Gamma$-flows on $D$
whenever $\lambda$ is equal to the order of $\Gamma$. 
Note that the definition of $F(G,\lambda)$ 
does not depend on the selection of 
$D$ and the additive Abelian group $\Gamma$ but on $G$ and 
the order of $\Gamma$. 
The function $F(G,\lambda)$ 
can also be obtained recursively by the following  rules
(see Tutte~\cite{tut}):
\begin{equation}\relabel{eq1-2}
F(G,\lambda)=
\left \{
\begin{array}{ll}
1, &\mbox{if }E=\emptyset;\\
0, &\mbox{if }G\mbox{ has a bridge};\\
F(G_1,\lambda)F(G_2,\lambda),
&\mbox{if }G=G_1\cup G_2;\\
(\lambda-1)F(G-e,\lambda), &\mbox{if }e\mbox{ is a loop};\\
F(G / e,\lambda)-F(G-e,\lambda),  \quad
&\mbox{if }e\mbox{ is not a loop nor a bridge},
\end{array}
\right.
\end{equation}
where $G/e$ and $G-e$
are the graphs obtained from 
$G$ by contracting $e$ and deleting $e$ respectively
and $G_1\cup G_2$ is the disjoint union of graphs 
$G_1$ and $G_2$.

A {\it block} of $G$ is 
a maximal subgraph of $G$ 
with the property that it is non-separable.
By (\ref{eq1-2}), the following result can be obtained.

\begin{lem}
\relabel{block-factor} 
If $G_1, G_2, \cdots, G_k$ are the components of $G$,
or $G_1, G_2, \cdots, G_k$ are the blocks of 
a connected graph $G$, then 
\beeq
F(G,\lambda)=\prod_{1\le i\le k}F(G_i,\lambda).
\eneq
\end{lem}

If $G$ is non-separable, $F(G,\lambda)$ can also be factorized 
when $G-e$ is separable for some edge $e$ 
or $G$ has a proper $3$-edge-cut $S$.
The results have been given in 
~\cite{jac3} (see \cite{dong2, jac2,jac4} also).

\begin{lem}[\cite{jac3}]
\relabel{v-edge}
Let $G$ be a bridgeless connected graph, $v$ be a vertex of $G$, 
$e = u_1u_2$ be an edge of $G$, and 
$H_1$ and $H_2$ be edge-disjoint subgraphs of 
$G$ such that $E(H_1) \cup E(H_2) = E(G - e)$, 
$V(H_1)\cap V(H_2) = \{v\}$, 
$V(H_1)\cup V(H_2) =V(G)$,
$u_1\in V (H_1)$ and $u_2\in V(H_2)$, as shown 
in Figure~\ref{f1}. 
Then 
\beeq
F(G, \lambda) =
\frac{F(G_1, \lambda)F(G_2, \lambda)}{\lambda -1}.
\eneq
where $G_i=H_i+vu_i$  for $i\in \{1, 2\}$.
\end{lem}

\begin{figure}[h!]
\centering 
\scalebox{0.9}{\input{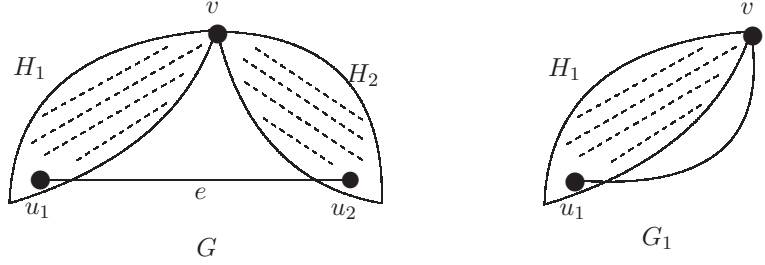}}
\caption{\relabel{f1} $G-e$ is separable.}
\end{figure}

If $G$ has an edge-cut $S$ with $2\le |S|\le 3$,
then $F(G,\lambda)$ also has a factorization \cite{jac3}. 

\begin{lem}[\cite{jac3}]
\relabel{2-edge}
Let $G$ be a bridgeless connected graph, 
$S$ be an edge-cut of $G$ 
and $H_1$ and $H_2$ be the sides of $S$, as shown in 
Figure~\ref{f2} when $|S|=3$. 
Let $G_i$ be obtained from $G$ by contracting 
$E(H_{3-i})$, for $i\in \{1, 2\}$. Then, for $2\le |S|\le 3$,
\beeq
F(G, \lambda) =
\frac{F(G_1, \lambda)F(G_2, \lambda)}
{(\lambda -1)_{|S|-1}},
\eneq
where $(x)_k$ is the polynomial $x(x-1)\cdots (x-k+1)$. 
\end{lem}

\begin{figure}[h!]
\centering 
\scalebox{0.9}{\input{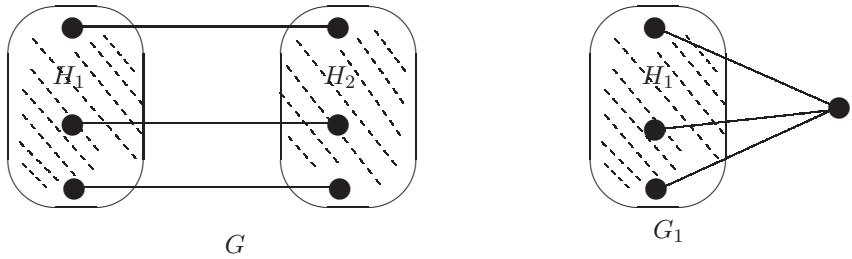}}
\caption{\relabel{f2} $G$ has a $3$-edge-cut.}
\end{figure}

It is not difficult to prove that for any 
bridgeless graph $G$, $F(G,\lambda)$ has no zero 
in $(-\infty, 1)$. 
But $1$ is a zero of $F(G,\lambda)$ whenever 
$G$ is not an empty graph. 
These conclusions can be obtained by equalities in (\ref{eq1-2}).
The next zero-free interval for flow polynomials
is $(1,32/27]$, due to Wakelin~\cite{wak}.

\begin{theo}[\cite{wak}]
\relabel{Wakelin}
Let $G=(V,E)$ be a bridgeless connected graph. Then
\begin{enumerate}
\item
$F(G, \lambda)$ is non-zero with sign 
$(-1)^{|E|-|V|+1}$ for $\lambda\in (-\infty, 1)$;
\item 
$F(G, \lambda)$ has a zero of multiplicity 
$1$ at $\lambda=1$ if $G$ is non-separable;
\item 
$F(G, \lambda)$ is non-zero for $\lambda\in (1, 32/27]$.
\proofend
\end{enumerate}
\end{theo}

For any integer $k\ge 0$, 
let $\xi_k$ be the supremum in $(1,2]$ 
such that $F(G,\lambda)$ is non-zero in the interval 
$(1,\xi_k)$ for all bridgeless graphs $G$ 
with at most $k$ vertices of degrees larger than 
$3$ (i.e., $|W(G)|\le k$). 
Clearly that $\xi_0\ge \xi_1\ge \xi_2\ge \cdots$. 
It is shown in \cite{dong2} that 
each $\xi_k$ can be determined by 
the flow roots of graphs from a finite set.

\begin{theo}[\cite{dong2}]
\relabel{dong-th1}
Each $\xi_k$ can be determined by 
the flow roots of graphs in a finite set $\Theta_k$, and 
$\xi_k=2$ for $k=0,1,2$,  
$\xi_3=1.430\cdots$, 
$\xi_4=1.361\cdots$ and $\xi_5=1.317\cdots$,
where the last three numbers 
are the real roots of $\lambda^3-5\lambda^2+10\lambda-7$,
$\lambda^3-4\lambda^2+8\lambda-6$ 
and $\lambda^3-6\lambda^2+13\lambda-9$
in $(1,2)$ respectively.
\end{theo}

\begin{cor}
\relabel{dong-cor}
$\xi_k>32/27$ for all $k\ge 0$.
\end{cor}

\proof 
By Theorem~\ref{dong-th1}, $\xi_k$ 
is determined by the flow roots of graphs from 
a finite set $\Theta_k$.
Thus $\xi_k$ is the flow root of some graph in $\Theta_k$.
By Theorem~\ref{Wakelin},
$\xi_k>32/27$.
\proofend

\resection{Graphs with real flow roots only\relabel{sec6}}

In this section, we assume that 
$G=(V,E)$ is a connected and bridgeless graph
and $r$, $\gamma, \alpha, k, \setr$ and $\omega$ are  
some invariants related to $G$ defined below:
\begin{enumerate}
\item  $r=|E|-|V|+1$;
\item  $\alpha=\sum\limits_{i\ge 3}(i-3)v_i$,
where $v_i$ is the number of vertices in $G$ of degree $i$;
\item  $\gamma$ is the number of $3$-edge-cuts of $G$;
\item $k=\sum\limits_{i\ge 4}v_i$;
\item $\setr$ is the multiset of real roots of $F(G,\lambda)$ 
in $(1,2)$; and 
\item $\omega=\sum\limits_{u\in \setr}(2-u)$.
\end{enumerate}
If we take another graph $H$, 
the above parameters related to $H$ are denoted by 
$r(H), \alpha(H), \gamma(H)$, $k(H), \setr(H)$
and $\omega(H)$ respectively.
It is straightforward to verify the following 
relations on these parameters. 

\begin{lem}\relabel{sect3-le1}
The following relations hold:
\begin{enumerate}
\item $\alpha=2|E|-3|V|$;
\item $k=W(G)=|V|-v_3$;
\item  $\gamma\ge v_3$, 
where the equality holds if and only if $G$ has no proper $3$-edge-cut;
\item $|V|=2r-2-\alpha$;
\item $|E|=3r-3-\alpha$;
\item $\omega+\sum_{u\in \setr}u=2|\setr|$.
\end{enumerate}
\end{lem}

It can be verified by (\ref{eq1-2}) that 
$F(G,\lambda)$ is a polynomial of order $r$.
Furthermore, if $G$ is $3$-edge connected, 
the coefficients of the three leading terms 
can be expressed in terms of $r, |E|$ and $\gamma$ 
(see \cite{kung}). 

\begin{lem}\relabel{le00}
If $G$ is $3$-edge connected, 
$F(G,\lambda)$ can be expressed as  
$\sum\limits_{0\le i\le r}b_i\lambda^i$, 
where $b_r=1, b_{r-1}=-|E|\mbox{ and }b_{r-2}={|E|\choose 2}-\gamma$.
\end{lem}
  
Recall that $L$ is the graph with one vertex 
and one loop and $\setgr$ is the family of  bridgeless 
graphs which have real flow roots only.
Obviously, we have the following conclusion on $r$. 

\begin{lem}
\relabel{le0}
$r\ge 1$.
Furthermore,  
if $G$ is a $3$-edge connected graph in $\setgr$, then 
$r=1$ if and only if $G$ is the graph $L$.
\end{lem}

From now on, we assume that $G$ is a $3$-edge connected 
graph in $\setgr$.  
By Lemma~\ref{le00}, 
we can get a lower bound for $\gamma$ 
in terms of $|E|$ and $r$.

\begin{lem}
 \relabel{le1}
Let $G=(V,E)$ be a $3$-edge connected graph in $\setgr$
with $|V|\ge 2$. 
Then
\begin{enumerate}
\item 
$\gamma\ge (|E|-r)(|E|-1)/(2r-2)$,
where the inequality is strict if $r-1$ does not divide $|E|-1$;

\item if $r\ge 3$ and $G$ is not an even graph\footnote{It is 
known that 
$G$ has a nowhere-zero $2$-flow if and only if every
vertex of $G$ has an even degree.}, then 
$\gamma\ge ((|E|-r)(|E|-4)+r-1)/(2r-4)$,
where the inequality is strict if $r-2$ does not divide $|E|-3$.
\end{enumerate}
\end{lem}

\proof  (i) 
It is known that any non-empty graph 
does not have nowhere-zero 1-flows, i.e., 
$F(G,\lambda)$ has a root $1$. 
Write 
\beeq\relabel{le-eq3}
F(G,\lambda)
=(\lambda -1)(\lambda^{r-1}-a_1\lambda^{r-2}
+a_2\lambda^{r-3}-\cdots).
\eneq
By Lemma~\ref{le00}, 
$a_1+1=|E|$ and $a_2+a_1={|E|\choose 2}-\gamma$.
So $\gamma={|E|\choose 2}-a_2-|E|+1$.
Since all roots of $F(G,\lambda)$ are real, 
applying Lemma 3.1 in \cite{kung} or 
the Newton Inequality \cite{har1978} to the coefficients of 
the three leading terms in the second factor of 
the right-hand side of (\ref{le-eq3}), 
we have 
\beeq\relabel{le-eq6}
a_2
\le {r-1\choose 2} 
\left (\frac{|E|-1}{r-1}\right )^2,
\eneq
where the inequality is strict if
$(|E|-1)/(r-1)$ is not a root of $F(G,\lambda)$.
Note that if $(|E|-1)/(r-1)$ is not an integer, 
it is not a root of $F(G,\lambda)$.
Thus 
$$
\gamma\ge {|E|\choose 2}-{r-1\choose 2} 
\left (\frac{|E|-1}{r-1}\right )^2-|E|+1
$$
and (i) follows.

(ii). It can be obtained similarly.  
As $G$ is not even, both $1$ and $2$ are flow roots of $G$.
Write 
\beeq\relabel{le-eq4}
F(G,\lambda)
=(\lambda -1)(\lambda -2)(\lambda^{r-2}-c_1\lambda^{r-3}
+c_2\lambda^{r-4}-\cdots).
\eneq
Applying the idea used in the proof of (i), we have 
$c_1=|E|-3$, ${|E|\choose 2}-\gamma=c_2+3c_1+2$ 
and 
$$
c_2\le {r-2\choose 2}\left (\frac{|E|-3}{r-2} \right )^2,
$$
where the inequality is strict whenever $\frac{|E|-3}{r-2}$ 
is not an integer. Thus
$$
\gamma\ge {|E|\choose 2}-{r-2\choose 2} 
\left (\frac{|E|-3}{r-2}\right )^2-3(|E|-3)-2
$$
and (ii) follows. 
\proofend

Recall that 
$\bar d(G)$ is the average value of degrees 
of all vertices $x\in W(G)$ in $G$,  i.e., 
$$
\bar d(G)=\frac{2|E|-3(|V|-k)}k.
$$

\begin{lem} \relabel{main-lem30}
Let $G=(V,E)$ be a $3$-edge connected graph in 
$\setgr$ with block number $b$. 
If $|V|\ge 3$ and $G$ does not contain any proper $3$-edge-cut, 
then  
\begin{enumerate}
\item  \relabel{lem30-no1} $r\ge 3$ and $|V|\ge 2k+1$;
\item \relabel{lem30-no2}
$|E|\ge 2|V|+2k-3+\frac{4(k-1)^2}{|V|-2k}$, 
where the inequality is strict if $r-2$ does not divide 
$|E|-3$; 
\item \relabel{lem30-no3} 
$|E|\ge |V|+8k-7$, 
where the inequality is strict 
if $|V|\ne 4k-2$;
\item \relabel{lem30-no4} $\bar d(G)\ge 
9+4\sqrt 2-2(\sqrt 2+1)^2/k
>14.656-11.656/k$;
\item \relabel{lem30-no5} $\omega\ge |E|-2|V|+2-b$,
where  
the inequality is strict if 
$F(G,\lambda)$ has some real roots in $(2,\infty)$;
\item \relabel{lem30-no6}  
$|E|\le \frac{\xi_k}{\xi_k-1} |V|+b-\frac{2}{\xi_k-1}
<b+\frac{32|V|-54}{5}$;
\item \relabel{cor1-no1} 
$k< \frac {11}{27}|V|+\frac{b+9}{54}$.
\end{enumerate}
\end{lem}

\proof 
(i) and (ii).
As $G$ is $3$-edge connected, 
$d(u)\ge 3$ for each vertex $u$ in $G$, implying that 
$$
r=|E|-|V|+1\ge \left \lceil 3|V|/2 \right \rceil -|V|+1
= \left \lceil |V|/2+1\right \rceil\ge 3.
$$ 
Since $G$ does not contain any proper $3$-edge-cut,
$\gamma=v_3=|V|-k$ by Lemma~\ref{sect3-le1}.
As $r=|E|-|V|+1$, by Lemma~\ref{le1} (i), we have 
\begin{equation}\relabel{lem300-eq1}
|V|-k\ge (|V|-1)(r+|V|-2)/(2(r-1)),
\end{equation}
which is equivalent to 
$$
(|V|-2k+1)r\ge |V|^2-|V|-2k+2.
$$
As $|V|\ge k$ and $|V|\ge 3$, we have $|V|^2-|V|-2k+2>0$.
Thus $|V|\ge 2k$, implying that 
$v_3=|V|-k\ge \frac {|V|}2>0$. 
Thus $G$ is not an even graph. 
As $r=|E|-|V|+1$, 
Lemma~\ref{le1} (ii) implies that  
\begin{equation}\relabel{lem300-eq2}
|V|-k\ge \frac{(|V|-1)(r+|V|-5)+r-1}{2(r-2)},
\end{equation}
where the inequality is strict if 
$r-2$ does not divide $|E|-3$. 
Observe that (\ref{lem300-eq2}) is equivalent to 
\begin{equation}\relabel{lem300-eq20}
r(|V|-2k)\ge |V|^2-4k-2|V|+4\ge 4(k-1)^2\ge 0.
\end{equation}
But $|V|^2-4k-2|V|+4=0$ implies that $k=1$ and 
$|V|=2k=2$, 
a contradiction. Thus $|V|\ge 2k+1$ and (i) holds.
The above inequality (\ref{lem300-eq20})
also implies that 
\begin{equation}\relabel{lem300-eq4}
r\ge \frac{|V|^2-4k-2|V|+4}{|V|-2k}
=|V|+2k-2+\frac{4(k-1)^2}{|V|-2k}.
\end{equation} 
If $r-2$ does not divide $|E|-3$,
then, by Lemma~\ref{le1} (ii),
the inequalities in  (\ref{lem300-eq2}) and 
(\ref{lem300-eq4}) are strict.
As $r=|E|-|V|+1$, the above inequality (\ref{lem300-eq4})
implies (ii) directly.

(iii).
By (\ref{lem300-eq4}), 
\begin{equation}\relabel{lem300-eq30}
r\ge 4k-2+|V|-2k+\frac{4(k-1)^2}{|V|-2k}
\ge 4k-2+2\times \sqrt {4(k-1)^2}=8k-6,
\end{equation}
where the last inequality is strict if and only if  
$|V|\ne 4k-2$.
As $r=|E|-|V|+1$, (iii) follows. 

(iv).  By (ii) and the definition of $\bar d(G)$, 
$$
4|V|+4k-6+\frac{8(k-1)^2}{|V|-2k}\le 2|E|=3(|V|-k)+k\bar d(G),
$$
implying that 
$$
\bar d(G)\ge 7+\frac{|V|-6+\frac{8(k-1)^2}{|V|-2k}}{k}
=7+\frac{2k-6+|V|-2k+\frac{8(k-1)^2}{|V|-2k}}{k}
$$
$$
\ge 7+\frac{2k-6+4\sqrt 2 (k-1)}{k}
=9+4\sqrt 2-2(\sqrt 2+1)^2/k
>14.656-11.656/k.
$$

(v).
Let $t=|\setr(G)|$, i.e., $t$ is the number of real roots 
of $F(G,\lambda)$ in the interval $(1,2)$.
Thus $t$ is the sum of the multiplicities of 
all flow roots of $G$ in $(1,2)$. 

By Theorem~\ref{Wakelin},
one root of $F(G,\lambda)$ is $1$ with multiplicity $b$, 
exactly $t$ of its roots are in $(1,2)$
and $(r-t-b)$ of its roots are at least $2$.
As $|E|$ is the sum of all flow roots of $G$, 
we have 
\beeq\relabel{main-lem1-eq1}
|E|\ge b+\sum_{u\in \setr}u+2(r-b-t)
=b+2t-\omega+2(r-b-t)
=2r-b-\omega,
\eneq
implying that $\omega\ge |E|-2|V|+1-b$ as $r=|E|-|V|+1$, 
where the inequality is strict if 
$F(G,\lambda)$ has some real roots in $(2,\infty)$.

(vi). 
Since $|V|\ge 2k+1$ by \ref{lem30-no1},
$G$ has some vertices of degree $3$ and thus 
$2$ is a root of $F(G,\lambda)$.
By Lemma~\ref{block-factor} and Theorem~\ref{Wakelin} (ii), 
$F(G,\lambda)$ has a root of multiplicity $b$ at $\lambda=1$. 
Thus $|\setr|\le r-1-b$ and 
\beeq\relabel{main-lem1-cor00-eq40}
|E|-|V|-b=r-b-1\ge |\setr|.
\eneq 
On the other hand, 
\beeq \relabel{main-lem1-cor00-eq41}
|\setr|(2-\xi_k)\ge \omega\ge |E|-2|V|+2-b,
\eneq 
where the last inequality is from \ref{lem30-no5}.
So (\ref{main-lem1-cor00-eq40}) and 
(\ref{main-lem1-cor00-eq41}) imply that 
\beeq   \relabel{main-lem1-cor00-eq42}
(|E|-|V|-b)(2-\xi_k)\ge |E|-2|V|+2-b.
\eneq 
Then it follows that 
\beeq \relabel{main-lem1-cor00-eq43}
|E|\le \frac{(|V|+b)\xi_k-b-2}{\xi_k-1}
=\frac{\xi_k}{\xi_k-1}|V|+\frac{b\xi_k-b-2}{\xi_k-1}
<b+\frac{32|V|-54}{5},
\eneq
where the last inequality follows from the fact that 
$\xi_k> 32/27$ by Corollary~\ref{dong-cor}.

(vii).  
By \ref{lem30-no2} and \ref{lem30-no6}, we have 
$$
2|V|+2k-3+\frac{4(k-1)^2}{|V|-2k}< b+(32|V|-54)/5.
$$
Solving this inequality gives that 
$k<\frac {11}{27}|V|+\frac{b+9}{54}$. So (vii) holds. 
\proofend

We are now going to establish 
the following important result. 
Recall that $\setgr_0$ is the family of 
non-separable and $3$-edge connected graphs 
$G$ in $\setgr$ 
such that $G$ does not contain 
any proper $3$-edge-cut
and $G-e$ is non-separable for each edge $e$ in $G$.

\begin{theo} \relabel{main-lem1}
Let $G=(V,E)\in \setgr_0$. 
If $G\not \in \{L, Z_3,K_4\}$, 
then  
\begin{enumerate}

\item \relabel{cor1-no01} 
$k\ge 3$;

\item \relabel{cor1-no2} 
$|\setr|\ge \nroot{k}\ge 9$,
where $\mu(x)$ is the function defined in 
Theorem~\ref{main-th} (ii).
\end{enumerate}
\end{theo}

\proof  
Suppose  that  $G\not \in \{L, Z_3,K_4\}$.
Clearly, $|V|\ge 2$, as $|V|=1$ and $G\in \setgr_0$ 
imply that $G=L$. 
It is easy to verify that all flow roots of 
$Z_s$ are real if and only if $s=3$.
As $G\ne Z_3$,  $|V|\ne 2$.
Hence $|V|\ge 3$.

(i)  
We first prove that $k\ne 0$.
Suppose that $k=0$.
Then $G$ is a cubic graph and so
$|E|=\frac 32 |V|$, implying that $|V|$ is even. 
Thus $|V|\ge 4$. 
By Lemma~\ref{main-lem30} (ii),
$|E|\ge 2|V|-2$. 
Thus $3|V|/2\ge 2|V|-2$, implying that $|V|\le 4$.
Thus  $|V|=4$. 
Since $G$ is cubic and $3$-edge connected,  
it is not difficult to verify that $G\cong K_4$, 
contradicting to the assumption.

Now suppose that  $1\le k\le 2$. 
By Theorem~\ref{dong-th1},
we have $\setr=\emptyset$ 
and so $\omega=0$. 
By Lemma~\ref{main-lem30} \ref{lem30-no5}, 
$|E|\le 2|V|-2+b=2|V|-1$ as $b=1$ (i.e., $G$ is non-separable).
By Lemma~\ref{main-lem30} \ref{lem30-no2},  
$$
|E|\ge 
\left \{ 
\begin{array}{ll}
2|V|-1, \quad &\mbox{if }k=1; \\
2|V|+2, &\mbox{if }k=2.
\end{array}
\right.
$$ 
Thus $k=1$ and 
$2|V|-1\le |E|\le 2|V|-1$,
implying that $|E|=2|V|-1$.
As $k=1$, $|E|=2|V|-1$ and 
$d(v)\ge 3$ for each vertex $v$ in $G$, 
it can be verified that 
$G$ has a vertex $u$ of degree $|V|+1$ 
and $|V|-1$ vertices of degree $3$.
So $G-u$ is a graph of order $|V|-1\ge 2$ and size $|V|-2$.
Since $G$ is non-separable, $G-u$ is connected.
Thus $G-u$ is a tree of order at least $2$, 
implying that $G-e$ is separable for each edge 
$e$ in $G-u$, contradicting the given condition
that $G\in \setgr_0$. 

Thus (i) holds.

(ii). By the definitions of $\xi_k$, $\omega$ and $\setr$, we have $\omega\le (2-\xi_k)|\setr|$.
As $k\ge 3$ and $G$ is non-separable, 
\ref{lem30-no5} and \ref{lem30-no2} in Lemma~\ref{main-lem30}
imply that $\omega\ge 2k-1$.
Thus  $|\setr| \ge (2k-1)/(2-\xi_k)$ 
by the definition of $\omega$.
By Theorem~\ref{dong-th1},
 $\xi_3\ge 1.430$, $\xi_4\ge 1.361$ and 
$\xi_5\ge 1.317$.
By Corollary~\ref{dong-cor}, $\xi_k> 32/27$ for all $k\ge 0$.
As $k\ge 3$ by (i), 
it is trivial to verify the result in \ref{cor1-no2}.
\proofend

Note that for $k=3,4,5,6,7,8,9,10$, 
the values of the function $\nroot{k}$ are 
respectively $9, 11,14,14, 16, 19, 21, 24$.

By Theorem~\ref{main-lem1} \ref{cor1-no2}, the following 
conclusion is obtained. 

\begin{cor}
\relabel{main-lem1-cor0}
Let $G=(V,E)\in \setgr_0$.
Then all flow roots of $G$ are integers if and only if 
$G\in \{L, Z_3,K_4\}$.
\end{cor}

Assume that $G=(V,E)\in \setgr_0$.
If some flow root of $G$ is not in the set $\{1,2,3\}$, 
then Theorem~\ref{main-lem1}~\ref{cor1-no01} and 
Lemma~\ref{main-lem30} \ref{lem30-no3} imply that 
$|E|\ge |V|+17$, 
and Theorem~\ref{main-lem1}~\ref{cor1-no2}
implies that $|\setr|\ge 9$.
In fact, these conclusions  still hold even if the condition
``$G\in \setgr_0$" is replaced by 
``$G\in \setgr$".

\begin{theo} \relabel{main-lem1-cor}
Let $G=(V,E)\in \setgr$. 
If some flow root of $G$ is not in the set $\{1,2,3\}$,
then $|E|\ge |V|+17$ and 
$|\setr|\ge 9$.
\end{theo}

\proof 
Let $\setz$ be the set of graphs in $\setgr$ 
which contain flow roots not in the set $\{1,2,3\}$.
Suppose that the result fails and 
$G$ is a graph in $\setz$ with the minimum value of 
$|E(G)|$ such that $|E|<|V|+17$ or $|\setr|<9$.
We first prove the following claims.

\inclaim $G$ is non-separable.

Suppose that $G$ is separable. 
By Lemma~\ref{block-factor}, 
some block $B$ of $G$ is contained in $\setz$.
By the minimality of $|E(G)|$, 
$|E(B)|\ge |V(B)|+17$ and $\setr(B)\ge 9$ hold.
As $G$ is bridgeless,
$|E(B')|\ge |V(B')|$ holds for each block of $G$.
Thus $|E(G)|\ge |V(G)|+17$ also holds.
By Lemma~\ref{block-factor} again,
$\setr(B)\ge 9$ implies that  $\setr(G)\ge 9$,
a contradiction. 
Thus  this claim holds. 

\inclaim $|V(G)|\ge 3$.

It is easy to verify that 
for any non-separable graph $H$ of order at most $2$, 
if all flow roots of $H$ are real, 
then each flow root of $G$ is in $\{1,2,3\}$.
As $G\in \setz$,  this claim holds.

\inclaim 
$G$ is $3$-edge connected.

Assume that $e$ is an edge 
contained in a $2$-edge-cut of $G$.  
By (\ref{eq1-2}), 
$F(G,\lambda)=F(G/e,\lambda)$.
Thus $\setr(G)=\setr(G/e)$, and 
$G\in \setz$ implies that  $G/e\in \setz$.
Also note that 
$|E(G)|-|V(G)|=|E(G/e)|-|V(G/e)|$,
implying that 
$G/e$ is also a counter-example to the result,
contradicting the assumption of $G$. 
Hence Claim~\thecountclaim\ holds.

\inclaim $G$ does not have any proper $3$-edge cut.

Suppose that $S$ is  a proper 
$3$-edge-cut of $G$, as shown in Figure~\ref{f2}. 
By Lemma~\ref{2-edge},
\beeq\relabel{main-lem1-cor-eq1}
F(G,\lambda)=\frac{F(G_1,\lambda)F(G_2,\lambda)}
{(\lambda -1)(\lambda-3)},
\eneq
where $G_1$ and $G_2$ are the 
graphs stated in  Lemma~\ref{2-edge}.
By (\ref{main-lem1-cor-eq1}), 
$G\in \setz$ implies that $G_i\in \setz$ for some $i$.
Say $i=1$. 
By the minimality of $|E(G)|$, 
$|E(G_1)|\ge |V(G_1)|+17$ and $|\setr(G_1)|\ge 9$ hold.
By (\ref{main-lem1-cor-eq1}) again, 
$|\setr(G_1)|\ge 9$ implies that $|\setr(G)|\ge 9$. 
As $G$ is bridgeless,  
it is not difficult to verify that 
$|E(G_1)|-|V(G_1)|\le |E(G)|-|V(G)|$.
Thus $|E(G_1)|\ge |V(G_1)|+17$  implies that 
$|E(G)|\ge |V(G)|+17$, a contradiction. 
Hence Claim~\thecountclaim\ holds. 

\inclaim $G-e$ is non-separable for each edge $e$ in $G$.

Suppose that $G-e$ is separable 
for some edge $e=u_1u_2$ as shown in Figure~\ref{f1}.  
By Lemma~\ref{v-edge}, 
\beeq\relabel{main-lem1-cor-eq2}
F(G,\lambda)=\frac{F(G_1,\lambda)F(G_2,\lambda)}
{\lambda -1},
\eneq
where $G_1$ and $G_2$ are the 
graphs stated in  Lemma~\ref{v-edge}. 
Then this claim can be proved similarly 
as the previous claim. 

By the above claims, we have $G\in \setgr_0\cap \setz$.
But, by Theorem~\ref{main-lem1} 
and Lemma~\ref{main-lem30}~\ref{lem30-no3}, we have
$k=|W(G)|\ge 3$,
$|E|\ge |V|+8k-7\ge |V|+17$
and $\setr(G)\ge 9$, contradicting the assumption of $G$. 

Hence the result holds.
\proofend

By Theorems~\ref{main-lem1-cor} and~\ref{kungroyle-theo},
we have the following result on plane graphs 
which have real chromatic roots only.

\begin{cor}
\relabel{cor3-th1}
Let $H$ be a 
connected planar graph of order $n$ and size $m$.
Assume that $H$ has real chromatic roots only 
but $H$ is not a chordal graph.
Then $n\ge 19$ and $H$ has at least $9$ chromatic roots in $(1,2)$
(counting multiplicity for each root). 
Furthermore, if every vertex-cut of $H$ 
does not induce a clique of $H$, 
then  $32n/27-5/9< m\le 2n-8$.
\end{cor}

\proof Assume that $H$ is a connected plane graph 
and $H^*$ is its dual. 
By the equality $P(H,\lambda)=\lambda F(H^{* },\lambda)$
due to Tutte~\cite{tut3},
the given conditions 
implies that $H^*$ has real flow roots only.
As $H$ is not chordal,
$P(H,\lambda)$ has non-integral roots
by the result in \cite{dong4} 
that planar graphs with integral
chromatic roots are chordal.
Thus $H^*$ has real flow roots only 
but also contains non-integral flow roots. 
By Theorem~\ref{main-lem1-cor}, 
$H^{* }$ has at least $9$ flow roots in $(1,2)$,
implying that 
$H$ has at least $9$ chromatic roots in $(1,2)$.
Notice that  $H^*$ has $m$ edges and $n$ faces. 
By Euler's polyhedron formula, the order of $H^*$ is  
\begin{equation}\relabel{sect3-eq5}
|V(H^*)|=m-n+2.
\end{equation}
By Theorem~\ref{main-lem1-cor} again, 
$|E(H^{* })|\ge |V(H^{* })|+17$, implying that 
$m\ge (m-n+2)+17$, i.e.,  $n\ge 19$. 

Now assume that every vertex-cut-set of $H$ does not 
induce a clique of $H$.
Then, it is trivial to verify that 
$H^*$ is a graph contained in $\setgr_0$.  
By Theorem~\ref{main-lem1}~\ref{cor1-no01}, $k(H^*)\ge 3$. 
Then, by  
Lemma~\ref{main-lem30} \ref{lem30-no2} and \ref{lem30-no6}, 
we have
\begin{equation}\relabel{sect3-eq06} 
2|V(H^*)|+4\le |E(H^*)|<(32|V(H^*)|-49)/5.
\end{equation}
By (\ref{sect3-eq5}), 
\begin{equation}\relabel{sect3-eq6}
2(m-n+2)+4\le m<(32(m-n+2)-49)/5.
\end{equation}
Thus $32n/27-5/9< m\le 2n-8$. 
\proofend

We end  this article with the following remark.

\noindent {\bf Remark}: 
By Lemmas~\ref{block-factor}, \ref{v-edge} and~\ref{2-edge},
the study of Problem~\ref{prob1} can be restricted to 
those graphs  in the family $\setgr_0$.
Thus, by Theorem~\ref{main-lem1}, 
there exist graphs asked in Problem~\ref{prob1} 
if and only if $\setgr_0-\{L, Z_3, K_4\}\ne \emptyset$.
By Theorem~\ref{main-lem1} again, 
for any $G\in \setgr_0-\{L, Z_3, K_4\}$, 
$G$ contains at least 
at least $\nroot{k}\ge 9$ 
flow roots in the interval $(1,2)$, where $k=|W(G)|\ge 3$.
However, as I know, 
no much research is conducted on 
counting the number of real flow roots of a graph 
in the interval $(1,2)$, 
except some study which confirms 
certain families of graphs having 
no real flow roots in the interval $(1,2)$
(see \cite{dong1, dong2, jac2,jac3,jac4}). 
It is unknown if there exists a graph $H$ 
with at least $\nroot{|W(H)|}$ flow roots in $(1,2)$.

\begin{prob}\relabel{prob2}
Is there a graph $H$ with $|W(H)|=k\ge 3$
and  
at least $\nroot{k}$ flow roots in $(1,2)$?
\end{prob}

\vspace{0.5 cm}

\subsection*{Acknowledgements}

The author wishes to thank the referees 
for their very helpful comments.

\end{document}